\documentclass{amsart}

\newtheorem{theorem}{Theorem}
\newcommand{\bt}{\begin{theorem}}
\newcommand{\et}{\end{theorem}}
\newtheorem{lemma}{Lemma}
\newcommand{\bl}{\begin{lemma}}
\newcommand{\el}{\end{lemma}}
\newtheorem{corollary}{Corollary}
\newcommand{\bc}{\begin{corollary}}
\newcommand{\ec}{\end{corollary}}
\newcommand{\beq}{\begin{equation}}
\newcommand{\eeq}{\end{equation}}
\newcommand{\benum}{\begin{enumerate}}
\newcommand{\eenum}{\end{enumerate}}

\newcommand{\R}{\ensuremath{\mathbf R}}

\newcommand{\bsmallmat}{\left(\begin{smallmatrix}}
\newcommand{\esmallmat}{\end{smallmatrix}\right)}

\newcommand{\bmat}{\left(\begin{matrix}}
\newcommand{\emat}{\end{matrix}\right)}

\newtheorem{problem}{Problem}
\newcommand{\bprob}{\begin{problem}}
\newcommand{\eprob}{\end{problem}}

\DeclareMathOperator{\sign}{\text{sign}}

\usepackage{amssymb,latexsym}

\title[Poincar\'e's Positivstellensatz]{Descartes's ``Rule of Signs'' and Poincar\'e's Positivstellensatz}

\author{Melvyn B. Nathanson}

\address{Department of Mathematics\\Lehman College (CUNY)\\Bronx, NY 10468}

\email{melvyn.nathanson@lehman.cuny.edu}



\date{\today}

\begin{document}

\begin{abstract} 
This is an exposition of Poincar\'e's 1883 paper, ``Sur les \'equations alg\'ebriques,'' 
which gives an important refinement of Descartes's rule of signs and was a precursor 
of P\' olya's Positivstellensatz. 

 \end{abstract}

\maketitle

\section{A refinement of Descartes rule of signs}


In 1637, Ren\' e Descartes~\cite{desc37} wrote,
\begin{quotation}
An equation can have as many true roots as it contains changes of sign, from $+$ to $-$ or 
from $-$ to $+$; and as many false roots as the number of times two $+$ signs or 
two $-$ signs are found in succession.
\end{quotation}
In the 17th century, a ``true root'' was a positive root and a ``false root'' was a negative root.  
Descartes's ``rule of signs'' is the most famous result in the theory of equations.  
Descartes did not prove the theorem.  In 1828, Gauss~\cite{gaus28,nath22-d} gave 
a beautiful elementary proof.  

In this paper we consider only polynomials with real coefficients.  
In more modern language, Descartes's theorem  
states that the number $Z(F)$ of positive roots (counting multiplicity) of a polynomial $F(x)$ 
 is the number $V(F)$ of sign variations in the sequence of coefficients 
of the polynomial minus a nonnegative even integer $\nu(F)$: 
\[ 
 Z(F) = V(F) - \nu(F).
\]
The number of sign variations in a sequence $(a_0, a_1, a_2, \ldots, a_n)$ 
of real numbers is the number of pairs $(j,k)$ with $0 \leq j < k \leq n$ 
such that $a_ja_k < 0$ and $a_i = 0$ if $j < i < k$.  
The number of sign variations of a polynomial $F(x) = \sum_{i=0}^n a_ix^i$ 
is the number of sign variations in the sequence $(a_0, a_1, a_2, \ldots, a_n)$ 
of coefficients of $F$. 
The sign variation function $V(F)$ gives an upper bound for the number of positive 
roots of $F(x)$, but only if $V(F)$ is odd does it imply that $F(x)$ has a positive root.  
The polynomial may or may not have positive root if $V(F)$ is even. 
For example, if $F(x) = x^2 -3x+2$, then $V(F) = Z(F) = 2$ and $\nu(F) = 0$.  
If $F(x) = x^2 -3x+5$, then $V(F) = \nu(F) = 2$ and $Z(F) =  0$.  
 It has been a nagging open problem in the theory of equations 
to understand the even integer $\nu(F)$.  
For example, virtual roots of polynomials~\cite{cost05,gonz98}  have been introduced 
to interpret $\nu(F)$.

Henri Poincar\' e was evidently also bothered by $\nu(F)$.  In a beautiful paper 
published in 1883, he applied an elementary argument to ``remove'' $\nu(F)$. 
Poincar\' e uses only algebraic and trigonometric identities 
to prove the following results.  
Proofs of these identities are collected in Appendix~\ref{Poincare:appendix}.

A polynomial has positive coefficients  
if all of its nonzero coefficients are positive.   
Sums and products of polynomials with positive coefficients 
are polynomials with positive coefficients. 
If $G(x)$ is a polynomial with positive coefficients, 
then $V(G) = 0$ and $G(x)$ has no positive root,

Poincar\' e~\cite{poin83} proved the following.

\bt                                           \label{Poincare:theorem:no-p}
A monic polynomial $F(x)$ is positive for all positive $x$ 
if and only if there exists a polynomial $G(x)$ with positive coefficients 
such that the product polynomial $F(x)G(x)$ has positive coefficients.  
\et

Equivalently, $F(x) > 0$ for all $x > 0$ if and only if 
there exists a nonzero polynomial $G(x)$ with positive coefficients 
such that $V(FG) = \nu(FG) =0$.   

In real algebraic geometry, a Positivstellensatz is a theorem 
that certifies that a polynomial is positive on some subset  of its domain.  
Theorem~\ref{Poincare:theorem:no-p} (Poincar\' e's Positivstellensatz) 
provides a certificate that a monic polynomial  is positive on 
$\Omega_1 = (0,\infty)$.  P\' olya~\cite{nath22-p,poly28} cites this result as a precursor 
of his Positivstellensatz for $n$-ary $m$-adic forms that are positive on 
the nonzero nonnegative orthant $\Omega_n$.

\bt                                                                 \label{Poincare:theorem:yes-p}
 Let $F(x)$ be a monic polynomial of degree $n$ 
with exactly $p$ positive roots  (counting multiplicity).  
Thus, $V(F) = p +  \nu(F)$ for some nonnegative even integer $\nu(F)$. 
There exists a nonzero polynomial $K(x)$  
such that $V(FK) = Z(FK) =p$ and $\nu(FK) = 0$. 
\et

\section{Proofs} 

Let  $\Lambda$ be the set of real and complex roots of the polynomial $F(x) \in \R[x]$. 
For all $\lambda \in \Lambda$, let $\mu_{\lambda}$ be the multiplicity 
of $\lambda$ as a root of the polynomial $F(x)$. 

If $\lambda =  \alpha + i\gamma$ is a complex root of $F(x)$ with $\gamma > 0$, then 
$\overline{\lambda} = \alpha - i \gamma$ is also a complex root of $F(x)$.  
We consider the quadratic polynomial $f_{\lambda}(x) \in \R[x]$ defined by 
\begin{align*}
f_{\lambda}(x) & =  (x-\lambda)(x-\overline{\lambda}) \\
& = x^2 - 2\alpha x + \alpha^2 + \gamma^2 \\ 
& = x^2 - 2\alpha x + \beta^2 
\end{align*}
where $\beta = \sqrt{\alpha^2 + \gamma^2} > |\alpha|$. 

Consider the following four subsets of $\Lambda$: 
\begin{align*}
\Lambda_1 & = \{\alpha \in \Lambda: \alpha \leq 0\} \\
\Lambda_2 & = \{\alpha + \gamma i \in \Lambda: \alpha \leq 0 \text{ and } \gamma > 0\} \\
\Lambda_3 & = \{\alpha + \gamma i \in \Lambda: \alpha > 0 \text{ and } \gamma > 0\} \\
\Lambda_4 & = \{\alpha \in \Lambda:  \alpha > 0 \}. 
\end{align*} 
Note that  the number of positive roots of $F(x)$ (counting multiplicity) 
is $\sum_{\alpha \in \Lambda_4} \mu_{\alpha}$.  
Associated with these sets are the monic polynomials 
\begin{align*}
F_1(x) & = \prod_{\alpha \in \Lambda_1} (x- \alpha)^{\mu_{\alpha}} \\
F_2(x) & = \prod_{\lambda \in \Lambda_2} f_{\lambda}(x)^{\mu_{\lambda}}   
 = \prod_{\lambda  = \alpha + \gamma i \in \Lambda_2} (x^2 - 2\alpha x + \beta^2)^{\mu_{\lambda}} \\
F_3(x) & = \prod_{\lambda \in \Lambda_3} f_{\lambda}(x)^{\mu_{\lambda}}   
= \prod_{\lambda  = \alpha + \gamma i \in \Lambda_3} (x^2 - 2\alpha x + \beta^2)^{\mu_{\lambda}} \\
F_4(x) & = \prod_{\alpha \in \Lambda_4} (x- \alpha)^{\mu_{\alpha}} 
\end{align*} 
and 
\[
F(x) = F_1(x) F_2(x) F_3(x) F_4(x).  
\]
We have $\alpha \leq 0$ for all roots in $\Lambda_1 \cup \Lambda_2$, and so $F_1(x)$ and $F_2(x)$ 
are polynomials with positive coefficients.

Let $\lambda = \alpha + \gamma i \in \Lambda_3$.  
We have $\gamma > 0$ and so 
\[
0 < \alpha < \sqrt{\alpha^2 + \gamma^2} = \beta. 
\]
There is a unique number $\varphi$ such that  
\[
0 < \varphi < \pi/2 
\]
 and 
\[
\cos \varphi = \frac{\alpha}{\beta}. 
\]
There is a unique integer $n \geq 2$ such that 
\[
0 < \varphi  < 2\varphi < \cdots < n\varphi < \pi \leq (n+1)\varphi < 3\pi/2.
\]
This inequality implies that 
\[
\sin (n+1)\varphi \leq 0 < \sin k \varphi  \qquad \text{for all $k \in \{1,2,\ldots, n\}$ }
\]  
Define the polynomial 
\[
g_{\lambda}(x) = \beta^{n-1} \sin \varphi  +  \left( \beta^{n-2} \sin 2\varphi \right) x +  \left( \beta^{n-3}  \sin 3\varphi \right) x^2
+ \cdots  +  \left( \sin n\varphi \right) x^{n-1}. 
\]
This is a polynomial of degree $n-1$ with positive coefficients.
From Lemma~\ref{Poincare:lemma:trigIdentity} in  Appendix~\ref{Poincare:appendix}  
we have the trigonometric identity  
\[
f_{\lambda}(x) g_{\lambda}(x)= \beta^{n+1}  \sin \varphi  -  \left(\beta \sin(n+1)\varphi \right) x^n +  \left(  \sin n\varphi \right) x^{n+1}.   
\]
It follows that $f_{\lambda}(x)g_{\lambda}(x)$ has positive coefficients. 
Let
\[
G(x) = \prod_{\lambda \in \Lambda_3} g_{\lambda}(x)^{\mu_{\lambda}}.   
\]
The polynomial 
\[
F_3(x) G(x) = 
\left( \prod_{\lambda \in \Lambda_3} f_{\lambda}(x)^{\mu_{\lambda}} \right) 
 \left(\prod_{\lambda \in \Lambda_3} g_{\lambda}(x)^{\mu_{\lambda}}  \right) 
=  \prod_{\lambda \in \Lambda_3} 
\left( f_{\lambda}(x) g_{\lambda} (x)\right)^{\mu_{\lambda}} 
\]
is a product of polynomials with positive coefficients, and so 
$F_3(x) G(x)$ has positive coefficients and 
\[
L(x) = F_1(x) F_2(x) F_3(x) G(x)
\]
also has positive coefficients.    

The set $\Lambda_4$ is the set of positive roots of $F(x)$.
If $\Lambda_4 = \emptyset$,  that is, if $F(x)$ has no positive root 
and $F_4(x) = 1$,  
then $F(x) = F_1(x) F_2(x) F_3(x)$ 
and $L(x) = F(x)G(x)$ has positive coefficients.  
This proves Theorem~\ref{Poincare:theorem:no-p}.

Suppose that $\Lambda_4 \neq \emptyset$. 
Let $q-1$ be the degree of the polynomial $L(x)$. 
For $\alpha \in \Lambda_4$, define the polynomial 
\[
h_{\alpha}(x) = \frac{x^q - \alpha^q}{x-\alpha} = \sum_{i=0}^{q-1} \alpha^{q-1-i} x^i.
\]
and let 
\[
H(x) = \prod_{\alpha \in \Lambda_4} h_{\alpha}(x)^{\mu(\alpha)}. 
\]
The product polynomial  
\[
M(x) = F_4(x) H(x) = \prod_{\alpha \in \Lambda_4} \left( (x-\alpha) h_{\alpha}(x) \right)^{\mu(\alpha)} 
= \prod_{\alpha \in \Lambda_4} \left( x^q - \alpha^q \right)^{\mu(\alpha)}
\]
is a monic polynomial of degree $pq$ that is a sum of powers of $x^q$.  Moreover, 
\[
V(M) = \sum_{\alpha \in \Lambda_4} \mu(\alpha) = p.
\]
by Lemma~\ref{Poincare:lemma:alternate}. 

We have 
\[
F(x)G(x) H(x) = F_1(x) F_2(x) F_3(x) G(x) F_4(x) H(x) =   L(x) M(x).
\]
Setting $K(x) = F(x)G(x)$ and applying Lemma~\ref{Poincare:lemma:VLM} gives 
\[
V(FK) = V(FGH) = V(LM) = V(M) = p. 
\]
This completes the proof. 

\appendix              
\section{Proofs of the identities}   \label{Poincare:appendix}

We prove the following trigonometric identity.  

\bl                          \label{Poincare:lemma:trigIdentity} 
Let
\[
f(x) = x^2 -  \left( 2\beta \cos \varphi \right) x + \beta^2
\]
and 
\[
g(x) = \beta^{n-1} \sin \varphi  +  \left( \beta^{n-2} \sin 2\varphi \right) x +  \left( \beta^{n-3}  \sin 3\varphi \right) x^2
+ \cdots  +  \left( \sin n\varphi \right) x^{n-1}. 
\]
Then 
\[
f(x) g(x)= \beta^{n+1}  \sin \varphi  -  \left(\beta \sin(n+1)\varphi \right) x^n +  \left(  \sin n\varphi \right) x^{n+1}. 
\]
\el

\begin{proof}
The trigonometric identities 
\[
\sin(k+1)\varphi = \sin k\varphi \cos \varphi + \cos k \varphi\sin \varphi
\]
and
\[
\sin(k-1)\varphi = \sin k\varphi \cos \varphi - \cos k\varphi\sin \varphi
\]
imply 
\beq                    \label{Poincare:trig-1} 
\sin(k+1)\varphi + \sin(k-1)\varphi =  2\sin k\varphi \cos \varphi.
\eeq
Also,
\beq                    \label{Poincare:trig-2} 
 \sin 2\varphi = 2 \sin \varphi \cos \varphi
 \eeq
 and 
 \begin{align}                   \label{Poincare:trig-3} 
\sin(n+1)\varphi & =  \sin n\varphi \cos \varphi + \cos n\varphi \sin \varphi  \nonumber \\ 
& =  2 \sin n\varphi \cos \varphi - \left(  \sin n\varphi \cos \varphi - \cos n\varphi \sin \varphi  \right)  \\
& =  2 \sin n\varphi \cos \varphi -  \sin (n-1) \varphi \nonumber 
 \end{align} 
Then 
\begin{align*}
f(x) g(x)
=  & \left(x^2 -  \left( 2\beta \cos \varphi \right) x + \beta^2 \right) 
 \left( \sum_{k=0}^{n-1} x^k \beta^{n-1-k}  \sin (k+1)\varphi \right)  \\ 
= & \  \beta^{n+1}  \sin \varphi  + \beta^n\left(  \sin 2\varphi -2 \sin \varphi \cos \varphi  \right)x \\
&  +  \ \sum_{k=2}^{n-1} \beta^{n+1-k} \left( \sin (k+1)\varphi - 2 \sin k \varphi \cos \varphi +  \sin (k-1)\varphi  \right) x^k \\
& - \ x^n \left( 2\sin n \varphi \cos \varphi  - \sin((n-1)\varphi \right) + x^{n+1}  \sin n\varphi. 
\end{align*}
Applying identities~\eqref{Poincare:trig-1},~\eqref{Poincare:trig-2}, and~\eqref{Poincare:trig-3}, 
we obtain 
\[
f(x) g(x)= \beta^{n+1}  \sin \varphi  - x^n \beta \sin(n+1)\varphi + x^{n+1}  \sin n\varphi. 
\]
This completes the proof. 
\end{proof}

\bl                                   \label{Poincare:lemma:alternate}
If $a_i > 0$ for all $i \in \{1,2,\ldots, n\}$, then 
\[
V\left( \prod_{i=1}^m (x-a_i) \right)  = m. 
\]
\el

\begin{proof}
Because $a_i > 0$ for all $i \in \{1,2,\ldots, n\}$, we obtain 
\begin{align*}
V\left( \prod_{i=1}^m (x-a_i) \right) & = V\left( \prod_{i=1}^m (x-\sign(a_i) ) \right) \\
&  = V\left(  (x-1)^m \right) = V\left( \sum_{j=0}^m (-1)^{m-j}  \binom{m}{j} x^j \right) \\
& = V\left( (-1)^m, (-1)^{m-1}, \ldots, 1,-1,1\right)  = m.
\end{align*}
This completes the proof. 
\end{proof}

\bl                                            \label{Poincare:lemma:VLM}
Let $M(x) = \sum_{i=0}^p r_ix^{qi}$ be a monic polynomial of degree $pq$ in which only powers of $x^q$ occur with nonzero coefficients.  
Let $L(x) = \sum_{j=0}^{q-1} s_j x^j$ be a polynomial of degree $q-1$ with positive 
coefficients.  Then
\[
V(LM) = V(M).
\]
\el

\begin{proof}
The sequence of length $(p+1)q$  of coefficients of the polynomial $M(x)$  is 
\[
\left( r_0, \underbrace{0,0,\ldots, 0}_{\text{$q-1$ roots} }, r_q,  \underbrace{0,0,\ldots, 0}_{\text{$q-1$ roots} }, r_{2q},\ldots, r_{(p-1)q}, \underbrace{0,0,\ldots, 0}_{\text{$q-1$ roots} }, r_{pq}, \underbrace{0,0,\ldots, 0}_{\text{$q-1$ roots} } \right). 
\]
For all $i \in \{ 0,1,2,\ldots, p\}$, we see the subsequence of coefficients of length $q$
\[
\left( r_{iq}, \underbrace{0,0,\ldots, 0}_{\text{$q-1$ roots} } \right). 
\]
In the sequence of coefficients of the product polynomial $L(x)M(x)$,  
this subsequence is replaced by 
\[
\left(  r_{iq} s_0, r_{iq} s_1,  r_{iq} s_2,  \ldots,  r_{iq} s_{q-1} \right). 
\]
For all $j \in \{0,1,2,\ldots, q-1$ we have $s_j > 0$ and so $\sign(r_{iq}s_j) = \sign(r_{iq})$.  
It follows that 
\[
V\left(  r_{iq} s_0, r_{iq} s_1,  r_{iq} s_2,  \ldots,  r_{iq} s_{q-1} \right) = 0
\] 
and 
\[
V(LM) = V(M).
\]
This completes the proof. 
\end{proof}


\begin{thebibliography}{99} 
\bibitem{cost05}
M. Coste, T.  Lajous-Loaeza, H. Lombardi, and M.-F. Roy, 
Generalized {B}udan-{F}ourier theorem and virtual roots,
J. Complexity 21 (2005), 479--486.
     
\bibitem{desc37} 
R. Descartes, \emph{The Geometry of Ren\' e Descartes}, 
Dover Publications, New York, 1954, page 373. 

\bibitem{gaus28}
C. F. Gauss, Beweis eines algebraischen Lehrsatzes, J. reine angew. Math. 3 (1828), 5--8.

\bibitem{gonz98}
L. Gonzalez-Vega, H.  Lombardi, and L. Mah\'{e}, Virtual roots of real polynomials,
J. Pure Appl. Algebra 124 (1998), 147--166. 



\bibitem{nath22-d} 
M. B. Nathanson, Gauss's proof of Descartes's rule of signs, arXiv: 2205.04249.

\bibitem{nath22-p} 
M. B. Nathanson, P\' olya's Positivstellensatz, preprint. 

\bibitem{poin83}
H. Poincar\' e, Sur les  \'equations alg\'ebriques, Comptes Rendus 97 (1883), 1418--1419.


\bibitem{poly28}
G. P\' olya, \" Uber positive Darstellung von Polynomen, Vierteljschr. Naturforsch. 
Ges. Z\" urich 73 (1928), 141--145; reprinted in G. P\' olya, 
\emph{Collected Papers, Vol. II: Location of Zeros}, MIT Press, Cambridge, 1974, pages 309--313. 
\end{thebibliography}
\end{document}